\documentclass{bcp}

\def\LaTeX{L\kern -.36em\raise .3ex\hbox{\sc a}\kern -.15em T\kern -.1667em%
\lower .7ex\hbox{E}\kern -.125em X}

\usepackage{bbm,epic}
\usepackage{amsmath,amssymb}

\newcommand{\ts}{\hspace{0.5pt}}
\newcommand{\NN}{\mathbb{N}}
\newcommand{\RR}{\mathbb{R}}

\newcommand{\Nnull}{\NN_{\ts 0}}
\newcommand\Ll[1]{\ell^{#1}}
\newcommand\vc[1]{\boldsymbol{#1}}  
\newcommand{\bs}{\boldsymbol}
\newcommand{\M}{\mathfrak{M}}
\newcommand{\cM}{\mathcal{M}}
\newcommand{\cR}{\mathcal{R}}
\newcommand{\m}{\mathfrak{m}}
\newcommand\dd[1]{\,\mathrm{d}#1}
\newcommand{\eps}{\varepsilon}

\newcommand\es\,  
\newcommand\is\,  

\newtheorem{thm}{Theorem}
\newtheorem{prop}{Proposition}
\newtheorem{lemma}{Lemma}

\begin{document}

\mathclass{Primary 92D25, 92D10; Secondary 60B10, 46B50.}
\thanks{It is a pleasure to thank E.\ Baake and O.\ Redner
for their cooperation, and C.\ Bank for carefully reading the
manuscript.}
\abbrevauthors{Michael Baake}
\abbrevtitle{Repeat distributions}

\title{Repeat distributions
from unequal crossovers}

\medskip

\author{Michael Baake}
\address{Fakult\"at f\"ur Mathematik, Universit\"at Bielefeld, \\
Postfach 100131, 33501 Bielefeld, Germany\\
E-mail: mbaake@math.uni-bielefeld.de}

\maketitlebcp

\abstract{It is a well-known fact that genetic sequences may contain
  sections with repeated units, called repeats, that differ in
  length over a population, with a length distribution of geometric
  type.  A simple class of recombination models with single crossovers
  is analysed that result in equilibrium distributions of this type.
  Due to the nonlinear and infinite-dimensional nature of these
  models, their analysis requires some nontrivial tools from measure
  theory and functional analysis, which makes them interesting also
  from a mathematical point of view. In particular, they can be viewed
  as quadratic, hence nonlinear, analogues of Markov chains.  }

\section*{1. Introduction.}

Recombination is a by-product of (sexual) reproduction, which leads to
the mixing of parental genes by exchanging genes (or sequence parts)
between homologous chromosomes (or DNA strands). This is achieved
through an alignment of the corresponding sequences, along with
crossover events that lead to a reciprocal exchange of the induced
segments. In this process, an imperfect alignment may result in
sequences that differ in length from the parental ones; this is known
as \emph{unequal crossover} (UC). Imperfect alignment is facilitated
by the presence of repeated elements (as is observed within some rDNA
sequences, compare \cite{GL}), and is believed to be an important
driving mechanism for the evolution of the corresponding copy number
distribution.  The perhaps best studied case of repeated elements
concerns microsatellites, see \cite{CS} and references given there for
a summary. An important observation is that, within a population, the
copy numbers vary, and often (at least approximately) follow a
distribution of geometric type (meaning a geometric distribution or a
finite convolution product thereof), see \cite{CS,KDSA,BWZBS} and
references therein for some experimental examples and findings.

The microsatellites themselves may follow an evolutionary course
independent of each other and thus give rise to evolutionary
innovation. For a detailed discussion of these topics, see
\cite{CS,SA} and references therein; for a brief introduction to
molecular evolution, see also \cite{Bus}, or \cite{Bue,Wake} for a
thorough overview. In this paper, which is mainly based on previous
work by Redner \cite{RB,Oli}, we shall focus on the distribution of
the copy numbers only, and disregard further aspects of the possible
evolution of the repeated units themselves. We rather aim at analysing
some simple models in order to understand the observed copy number or
repeat distributions.  Moreover, we are primarily interested in models
that preserve the mean copy number, though our setting will be
adequate to accommodate also more general models. In view of possible
applications to systems where the copy number (slowly) changes with
time, it seems natural to set up a frame that can cope with such a
situation as well.

In the entire model class to be described below, one considers
individuals whose genetic sequences contain a section with repeated
units.  These may vary in number, $i\in\NN_{0} = \{0,1,2,3,\ldots\}$,
where $i=0$ is explicitly allowed and corresponds to no unit being
present (yet).  The composition of these sections (with respect to
mutations that might have occurred) and the rest of the sequence are
ignored here, as are details of the actual alignment process (e.g.,
whether partial loops of longer pieces are formed in order not to
disturb the alignment outside the repeat region), see also \cite{Bank}
for a first discussion of possible models in this direction.

\begin{figure}
\medskip
  \begin{center}
    \unitlength=1mm
     \begin{picture}(85,11)(0,0)
      \matrixput(0.625,9.5)(2.5,0){4}(0,0){1}{\line(1,0){1.25}}
      \matrixput(10,8)(5,0){9}(0,3){2}{\line(1,0){5}}
      \matrixput(10,8)(5,0){10}(0,3){1}{\line(0,1){3}}
      \matrixput(60,8)(5,0){3}(0,3){2}{\line(1,0){5}}
      \matrixput(60,8)(5,0){4}(0,3){1}{\line(0,1){3}}
      \matrixput(75.625,9.5)(2.5,0){4}(0,0){1}{\line(1,0){1.25}}
      \matrixput(0.625,1.5)(2.5,0){10}(0,0){1}{\line(1,0){1.25}}
      \matrixput(25,0)(5,0){6}(0,3){2}{\line(1,0){5}}
      \matrixput(25,0)(5,0){7}(0,3){1}{\line(0,1){3}}
      \matrixput(60,0)(5,0){2}(0,3){2}{\line(1,0){5}}
      \matrixput(60,0)(5,0){3}(0,3){1}{\line(0,1){3}}
      \matrixput(70.625,1.5)(2.5,0){6}(0,0){1}{\line(1,0){1.25}}
      \qbezier(55,9.5)(57.5,5.5)(59.9,1.5)
      \qbezier(55,1.5)(57.5,5.5)(59.9,9.5)
    \end{picture}
\medskip
\caption{Snapshot after an unequal crossover event as described in the
  text.  Rectangles denote the relevant blocks, while the dashed lines
  indicate possible extensions with other elements that are
  disregarded here.}
    \label{fig:ureco}
  \end{center}
\end{figure}
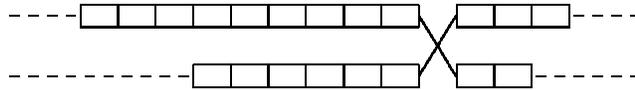

In the course of time, recombination events take place in which a
random pair of individuals is formed and their respective sections are
randomly aligned, possibly imperfectly with `overhangs'.  Then, both
sequences are cut at a common position between two building blocks and
their right (or left) fragments are interchanged.  This so-called
\emph{unequal crossover} is schematically depicted in
Figure~\ref{fig:ureco}.  Obviously, the total number of relevant units
is conserved in each event. While this is clearly a stochastic
process, it is nevertheless interesting to investigate its
deterministic limit, at least as a first step towards a better general
understanding of this model class.  To contribute to this first step,
and to summarise what has been done in this direction so far, is the
main aim of this contribution.

\section*{2. Description of the deterministic limit.}

As a first step for the analysis of crossover dynamics, we assume the
population size to be (effectively) infinite, i.e., large enough so
that random fluctuations may be neglected (finite populations will
briefly be mentioned later on).  We write $\cM (X)$ for the (finite)
measures on a space $X$, denote the restriction to positive measures
by a superscript $+$, and indicate a restriction to measures of total
variation $r$ by a corresponding subscript (see \cite{ReSi} or
\cite{Yos} for a short summary of the measure theory needed here).
Then, the distribution of the copy numbers over our population is
described by a probability measure (or vector) $\vc{p} \in
\cM_1^+(\Nnull)$, which we identify with an element $\vc{p} =
(p_k^{})_{k\in\Nnull}^{}$ in the appropriate subset of
$\Ll{1}(\Nnull)$.  Since we do not consider any genotype space other
than $\Nnull$ in this article, reference to it will be omitted in the
sequel, so we write $\ell^1$ instead of $\ell^1 (\Nnull)$ from now on.
These spaces are complete in the metric derived from the usual
$\Ll{1}$ norm, which is the same as the total variation norm here.
The metric is denoted by
\begin{equation} \label{urgpmetric}
  d(\vc{p},\vc{q}) \; = \;
  \|\vc{p} - \vc{q}\|^{}_1 \; = \;
   \sum_{k\ge0} |p_k^{} - q_k^{}| \,.
\end{equation}

Let us consider the above process (as well as various more general
ones) on the level of the induced dynamics on the probability measures
(i.e., in the infinite population limit mentioned above).  With the
notation just introduced, the dynamics can be described by means of
the \emph{recombinator}
\begin{equation} \label{reco}
  \cR(\vc{p})_i^{} \; := \; 
  \frac{1}{\|\vc{p}\|^{}_1}
    \sum_{j,k,\ell\ge0} \, T_{ij,k\ell}^{} \, p_k^{} \, p_\ell^{} \,.
\end{equation}
Here, $T_{ij,k\ell} \ge 0$ denotes the probability that a pair
$(k,\ell)$ turns into $(i,j)$, so, for normalisation, we require
\begin{equation} \label{sumt}
  \sum_{i,j\ge0} T_{ij,k\ell} \; = \; 1 \, ,
  \qquad\text{for all $k, \ell \in \Nnull$.}
\end{equation}
The factor $p_k^{} \, p_\ell^{}$ in \eqref{reco} describes the
probability that a pair $(k,\ell)$ is formed, i.e., we assume that two
individuals are chosen independently from the population.  We assume
further that, for all $i,j,k,\ell$, 
\begin{equation} \label{symmetry}
   T_{ij,k\ell} = T_{ji,k\ell} \; = \; T_{ij,\ell k} \, ,
\end{equation}
i.e., that $T_{ij,k\ell}$ is symmetric with respect to both index
pairs, which is reasonable and follows from the corresponding symmetry
of the underlying process, compare Figure~\ref{fig:ureco}.  Then, the
summation over $j$ in \eqref{reco} represents the breaking-up of the
pairs after the recombination event.  These two ingredients (symmetry
and summation) lead to the quadratic nature of the iteration process,
see below for more and \cite{L1,L2} for the appearance of similar
types of equations in a different class of biological models.

Condition \eqref{sumt} and the presence of the prefactor
$1/\|\vc{p}\|^{}_1$ in the defining Eq.~\eqref{reco} make $\cR$ norm
non-increasing, i.e., $\|\cR(\vc{x})\|^{}_1 \le \|\vc{x}\|^{}_1$, and
positive homogeneous of degree 1, i.e., $\cR(a\vc{x}) = |a|
\cR(\vc{x})$, for $\vc{x} \in \Ll{1}$ and $a \in \RR$.  Furthermore,
$\cR$ is a positive operator with $\|\cR(\vc{x})\|^{}_1 =
\|\vc{x}\|^{}_1$ for all positive elements $\vc{x} \in \Ll{1}$.  Thus,
it is guaranteed that $\cR$ maps $\cM_r^+$, the space of positive
measures of total variation $r$, into itself. This subspace is
complete in the topology induced by the norm $\|.\|_1^{}$, i.e., by
the metric $d$ from \eqref{urgpmetric}.  (For $r=1$, the prefactor on
the right hand side of \eqref{reco} is redundant, but improves
numerical stability of an iteration with the nonlinear mapping $\cR$.)

Given an initial configuration $\vc{p}_0^{} = \vc{p}(0)$,
the dynamics may be taken in discrete time steps, with subsequent
generations,
\begin{equation}
  \label{urdisctime}
  \vc{p}(t+1) \; = \; \cR(\vc{p}(t)) \,,
  \qquad t\in\Nnull \,.
\end{equation}
This iteration reflects the following: due to random mating, it is
sufficient to consider the dynamics at the level of the single
strands, which will be combined into pairs again randomly in the next
generation, according to the Hardy-Weinberg equilibrium \cite{Bue}.

Our treatment of this case will be set up in a way that also allows
for a generalisation of the results to the analogous process in
continuous time, where generations are overlapping,
\begin{equation}
  \label{urconttime}
  \tfrac{\dd{}}{\dd{t}}\, \vc{p}(t) \; = \; \varrho \, 
   (\cR-\mathbbm{1})(\vc{p}(t)) \,,
  \qquad t \ge 0 \,.
\end{equation}
This reflects what is called \emph{instant mixing}, i.e., the
instantaneous formation of pairs, their recombination and separation.
In other words, the actual duration of the diplophase (or ``paired
phase'') is neglected, which is an approximation that is justified
as long as recombination is rare on the time scale of an
individual life span.

Obviously, the (positive) parameter $\varrho$ in \eqref{urconttime}
only leads to a rescaling of the time $t$.  We therefore choose
$\varrho=1$ without loss of generality.  Furthermore, it is easily
verified that the fixed points of \eqref{urdisctime} are in one-to-one
correspondence with the equilibria of \eqref{urconttime}.  (In the
sequel, we use the term fixed point for both discrete and continuous
dynamics.)

As mentioned above, our main interest at present is in processes that
conserve the total copy number in each event, i.e., $T^{}_{ij,k\ell} >
0$ for $i+j = k+\ell$ only. More general scenarios are possible, and
also interesting, but already the concept of an equilibrium gets
rather involved, whence we do not go into further details here.
To\-geth\-er with the normalisation \eqref{sumt} and the symmetry
condition from above, this yields
\begin{equation} \label{urcons}
  \sum_{i,j\ge0} i \, T_{ij,k\ell}
  \; = \; \sum_{i,j\ge 0} \frac{i+j}{2} \, T_{ij,k\ell}
  \; = \; \sum_{i,j\ge 0} T_{ij,k\ell} \, \frac{k+\ell}{2}
  \; = \; \frac{k+\ell}{2} \,,
\end{equation}
the second equality of which is an alternative condition that
can replace the strict preservation of the copy number as follows.

\smallskip
\begin{lemma} \label{mean-preservation} Let $\cR$ be defined by\/
  $\eqref{reco}$, with $T_{ij,k\ell}\ge 0$ subject to the
  normalisation\/ $\eqref{sumt}$ and the symmetry conditions\/
  $\eqref{symmetry}$. If also the second equality in\/
  $\eqref{urcons}$ is satisfied, for all $k,\ell\in\Nnull$, the mean
  copy number in the population is preserved.
\end{lemma}
\Proof This is a simple calculation,
\[
  \sum_{i\ge0} i \, \cR(\vc{p})_i^{} \, = \, 
  \sum_{i,j,k,\ell\ge0} i \, T^{}_{ij,k\ell} \, p_k^{} \, p_\ell^{}
  \, = \, \sum_{k,\ell\ge0} \frac{k+\ell}{2} \, p_k^{} \, p_\ell^{}
  \, = \, \sum_{k\ge0} k \, p_k^{} \, ,
\]
which shows the claim, provided that the mean
$\mathfrak{m}:=\sum_{i} i \ts p_i$
is well-defined. \hfill \sq

{}From now on, we use the symbol $\mathfrak{m}$ for the mean,
in order not to confuse it with summation indices and the like.

\section*{3. Markov chains for comparison.}

Let us take a brief detour to look at the linear counterpart, a
countable state Markov chain, in the deterministic limit of the
forward equation for the time evolution of its probability
distribution.  To this end, consider again probability vectors
$\bs{p}$ on $\Nnull$ and define
\[
    \M (\bs{p})^{}_k \; := \; 
    \sum_{\ell=0}^{\infty} M^{}_{k \ell}\, p^{}_\ell \, ,
\]
for $k\in\Nnull$, where all $M_{k \ell}\ge 0$ together with
$\sum_{k=0}^{\infty} M_{k \ell} = 1$ for all $\ell\in\Nnull$.  This
also makes the above sums well-defined on all elements of $\ell^{1}$.
Note that the matrix $M = (M_{k\ell})_{k,\ell\in\Nnull}$ is transposed
in comparison with the standard convention for Markov chains \cite{S},
because we are using it here in a dynamical systems context, with
action of the matrix to the column vector on the right.  The time
evolution now either reads
\begin{equation} \label{markov-discrete}
   \bs{p}\ts (t+1) \; = \; \M (\bs{p} (t))
   \qquad \mbox{(in discrete time)}
\end{equation}
or
\begin{equation} \label{markov-cont}
   \tfrac{\dd}{\dd t}\bs{p} \; = \; (\M - \mathbbm{1}) (\bs{p})
   \qquad \mbox{(in continuous time),}
\end{equation}
where the rate constant is again assumed to be $1$, compare
the remark after Eq.~\eqref{urconttime}.

The iteration of the discrete version \eqref{markov-discrete} on
$\ell^{1}$ is well-defined, while uniqueness of the solution of the
initial value problem for the continuous time counterpart
\eqref{markov-cont} on the same space follows from its global
Lipschitz property,
\[
   \left\| \M (\bs{p}) - \M (\bs{q}) \right\|^{}_{1}
  \, = \, \left\| \M (\bs{p} - \bs{q}) \right\|^{}_{1}
  \, \le \, \sum_{k,\ell \ge 0} M^{}_{k \ell}\, \lvert p^{}_{\ell}
  - q^{}_{\ell}\rvert \, = \, \sum_{\ell \ge 0} \, \lvert
  p^{}_{\ell} - q^{}_{\ell}\rvert \, = \,
  \left\| \bs{p} - \bs{q} \right\|^{}_{1} \, ,
\]
which holds for all $\bs{p},\bs{q}\in\ell^1$. The properties of the
matrix $M$ guarantee that the positive cone as well as the simplex of
probability vectors are preserved in forward time.  Consequently, one
can consider \eqref{markov-discrete} and \eqref{markov-cont} as
dynamical systems on $\ell^1$. As the latter is a Banach space of
infinite dimension, the unit ball is no longer compact in the norm
topology, whence some extra care is needed for the results.

As before, fixed points of \eqref{markov-discrete} line up with
equilibria of \eqref{markov-cont}, so that we speak of fixed points in
both cases. Their existence is provided by Perron-Frobenius theory for
countable state Markov matrices, see \cite[Ch.~7.1]{K} or
\cite[Ch.~5]{S} for a detailed account. Irreducibility, aperiodicity
and primitivity are defined as in the finite-dimensional case without
difficulty. However, for meaningful results on eigenvalues and
eigenvectors, one additionally needs the concept of \emph{recurrence},
see \cite[p.~197 f.]{K} for a nice summary.

The Perron value $\lambda$ emerges from the radius of convergence,
$\rho$, of the power series $T(z)=\sum_{n\ge 0} (zM)^n$ via $\rho =
1/\lambda$.  Clearly, we have $\rho\ge 1$ for a Markov matrix.  If one
diagonal entry (and then any) of $T(z)$ diverges at $1$ (so that
$\rho=1$ in this case, compare \cite[Thm.~6.6]{S}), the countable
state Markov matrix $M$ is called \emph{recurrent}, where the
behaviour of the diagonal element $T(z)_{ii}$, as $z\to 1$, reflects
the expected number of recurrences to $i$, which is infinite in this
case.  Moreover, a unique normalised and strictly positive (right)
eigenvector $\vc{p}\in\cM^{+}_{1}$ exists with $\M(\vc{p})=\vc{p}$,
see \cite[Thm.~5.4]{S}. This probability vector has the meaning of the
unique equilibrium distribution and is the desired fixed point of the
dynamics.

Assume for a moment, in addition to the above conditions on $\M$, that
\begin{equation} \label{markov-mean-one}
    \sum_{i\ge 0} i M_{ij} \; = \; j \, , \quad \mbox{for all $j$}.
\end{equation}
As before, this is a sufficient condition for the mean to be preserved
under the dynamics, because one has
\[
   \sum_{i\ge 0}  i\, \M (\bs{p})_i \, = \, \sum_{i\ge 0} \sum_{j\ge 0}
   i M_{ij} p_j \, = \, \sum_{j\ge 0} \sum_{i\ge 0}
   i M_{ij} p_j \, = \, \sum_{j\ge 0} j\ts p_j \, ,
\]
with the interchange of summation being permissible due to absolute
convergence of the sums involved, provided that $\m = \sum_j j p_j$
exists. However, a condition of type \eqref{markov-mean-one} is
usually too restrictive for a linear system, wherefore we do not
impose it here. As we shall see, the mean copy number can be preserved
without it.

A probability vector $\bs{p}$ is called \emph{reversible} for $\M$
when the detailed balance equation
\begin{equation} \label{markov-rev}
    M^{}_{k \ell} \, p^{}_{\ell} \; = \; M^{}_{\ell k} \, p^{}_k
\end{equation}
holds for all $k,\ell\in\Nnull$. An important consequence is
that any reversible $\bs{p}$ is automatically a fixed point
of $\M$:
\[
   \M (\bs{p})^{}_k \, = \, 
   \sum_{\ell\ge 0} M^{}_{k \ell} \, p^{}_{\ell}
   \, = \, \sum_{\ell\ge 0} M^{}_{\ell k} \, p^{}_k \, = \, p^{}_k \, .
\]
Reversibility often provides a simpler way to actually calculate a
specific fixed point than the defining matrix eigenvalue equation.

Since the Perron-Frobenius eigenvalue $\lambda$ need not be isolated
in the spectrum of $M$, the convergence properties are more subtle
than in the finite-dimensional situation. Under certain extra
conditions (e.g., if $\lambda$ \emph{is} isolated), the time evolution
of an arbitrary initial condition converges exponentially fast towards
the fixed point. However, when the matrix $M$ is not only recurrent,
but also positive recurrent, one has at least convergence of the
discrete iteration, see \cite{S} for details. Here, positive
recurrence means that the expected time for a return to the state $i$
is \emph{finite}, which is clearly stronger than mere recurrence.
     
\smallskip
The standard geometric distribution with parameter $\alpha\in (0,1)$ 
is a discrete probability distribution on $\Nnull$, defined by the
probability vector $\vc{p}$ with
\begin{equation}  \label{geom-def}
   p_n \, := \, \alpha\ts (1-\alpha)^n \; , \quad 
   \mbox{for $n\in\NN_0$}\, .
\end{equation}
Clearly, $p_n > 0$ and $\sum_{n\ge 0} p_n = 1$, while $\m = \sum_{n\ge
  0} n \ts p_n = (1-\alpha)/\alpha$, so that $\alpha =
1/(\mathfrak{m}+1)$.  If we define the matrix $M = (M_{ij})_{i,j\ge
  0}$ by $M_{ij} = p_i$, one has
\[
   (M\vc{p})_i \, = \, \sum_j M_{ij} p_{j} \, = \,
   p_i \sum_{j} p_j \, = \, p_i \, ,
\]
so that $M\vc{p} = \vc{p}$. One clearly has $M^n=M$ for all
$n\in\NN$. Consequently, each entry of $T(z)$  is a geometric
series of the form $M_{ij} (1+z+z^2+\ldots)$, which thus
diverges at $z=1$. In particular, $M$ is (positive)
recurrent.

The matrix $M$ does not satisfy Eq.~\eqref{markov-mean-one}.
Nevertheless, the mean copy number is preserved in the following
sense. Let $\vc{a}$ be an arbitrary probability vector with mean $\m$,
and $\vc{p}$ the geometric distribution according to \eqref{geom-def}
with the same mean. With the corresponding matrix $M$, one then finds
\[
    \sum_{i,j} i\ts M_{ij}\ts a_j \, = \,
    \sum_{i,j} i\ts p_i \ts a_j \, = \,
    \sum_{i} i\ts p_i \,\sum_{j} a_j \, = \, \m\, ,
\]
which results in the mean preservation, provided one starts with an
initial condition $\vc{a}$ of mean $\m$. Otherwise, the iteration maps
$\vc{a}$ to an image of mean $\m$ in the first step, and preserves
$\m$ in all subsequent iterations.

{}Further eigenvectors of $M$ are given by $\vc{q}^{(\ell)}:=
\vc{e}^{}_{0} - \vc{e}^{}_{\ell}$ for $\ell\in\NN$, where $\vc{e}_i$
is the standard basis vector with $1$ in coordinate $i$ and $0$
otherwise. All these extra vectors belong to the eigenvalue $0$, which
is the only other eigenvalue of $M$. In fact, $M$ is diagonalisable,
and it is not difficult to see that an arbitrary vector $\vc{a}=
(a^{}_{0},a^{}_{1},a^{}_{2},\ldots) \in \ell^{1}$ can be written as a
convergent expansion, $\vc{a} = \beta\ts\vc{p} + \sum_{\ell\ge 1}
(\beta\ts p^{}_{\ell} - a^{}_{\ell})\ts\vc{q}^{(\ell)}$, where $\beta
= \sum_{i\ge 0} a_{i}$.  Consequently, the chosen eigenvectors of $M$
form a basis of $\ell^{1}$. If $U =
(\vc{p},\vc{q}^{(1)},\vc{q}^{(2)},\ldots)$ denotes the matrix that
columnwise consists of the eigenvectors of $M$, one has
\[
   M \, = \, U\ts \mbox{diag} (1,0,0,\ldots)\, U^{-1},
\]
which makes the relation $M^n=M$ for $n\in\NN$ particularly
transparent. Moreover, one sees that $M$ commutes with all matrices
$N$ of the form $N = U\! A\ts U^{-1}$ where $A$ has the block form
\[
    A \, = \, \begin{pmatrix} a & \vc{0}^t \\
              \vc{0} & A'  \end{pmatrix} 
\]
with an arbitrary matrix $A'$. Restricting $N$ so that $M+N$ is still
Markov, one can find multi-parameter families of Markov matrices that
share the given stationary geometric distribution $\vc{p}$. The same
stationary probability vector $\vc{p}$ can thus arise from many other
Markov chains as well.

Let us now return to the bilinear counterpart to see which
of these structural properties possess an analogue, and to
describe the setting of our later analysis.

\section*{4. General structure of the bilinear system.}

Consider the crossover dynamics as defined by \eqref{reco}.  Let us
begin by stating the following general fact.

\smallskip
\begin{prop}
  \label{recolipschitz}
  If the recombinator\/ $\cR$ of\/ $\eqref{reco}$ satisfies  the
  normalisation conditions\/ $\eqref{sumt}$, one has
  the global Lipschitz condition
  \begin{displaymath}
    \|\cR(\vc{x}) - \cR(\vc{y})\|^{}_1 \; \le \; C 
    \|\vc{x} - \vc{y}\|^{}_1 \, ,
  \end{displaymath}
  with constant\/ $C=3$ on\/ $\Ll{1}$, respectively\/
  $C=2$ if\/ $\vc{x}$, $\vc{y} \in \cM_r$.
\end{prop}
\Proof Let $\vc{x}$, $\vc{y} \in \Ll{1}$ be non-zero (otherwise the
  statement is trivial).  Then, one has
\begin{align*}
  \lefteqn{  \|\cR(\vc{x}) - \cR(\vc{y})\|^{}_1  \; = \;
    \sum_{i\ge0}\; \biggl|
      \sum_{j,k,\ell\ge0} T_{ij,k\ell}^{} \left(
        \frac{x_k^{}\,x_\ell^{}}{\|\vc{x}\|^{}_1} -
        \frac{y_k^{}\,y_\ell^{}}{\|\vc{y}\|^{}_1} \right) \biggr|} \\
    &\; \le \; \sum_{k,\ell\ge0} \left|
      \frac{x_k^{}\,x_\ell^{}}{\|\vc{x}\|^{}_1} -
      \frac{y_k^{}\,y_\ell^{}}{\|\vc{y}\|^{}_1} \right|
    \sum_{i,j\ge0} T_{ij,k\ell} 
    \; = \; \sum_{k,\ell\ge0} \left|
      \frac{x_k^{}\,x_\ell^{}}{\|\vc{x}\|^{}_1} -
      \frac{x_k^{}\,y_\ell^{}}{\|\vc{x}\|^{}_1} +
      \frac{x_k^{}\,y_\ell^{}}{\|\vc{x}\|^{}_1} -
      \frac{y_k^{}\,y_\ell^{}}{\|\vc{y}\|^{}_1} \right| \\
    &\; \le \; \sum_{k,\ell\ge0} \left(
      \frac{|x_k^{}|}{\|\vc{x}\|^{}_1} |x_\ell^{} - y_\ell^{}| +
      |y_\ell^{}| \left| \frac{x_k^{}}{\|\vc{x}\|^{}_1} - 
        \frac{y_k^{}}{\|\vc{y}\|^{}_1} \right|\, \right) 
    \; = \; \|\vc{x}-\vc{y}\|^{}_1 + 
    \frac{\bigl\|\|\vc{y}\|^{}_1 \vc{x} - 
     \|\vc{x}\|^{}_1 \vc{y} \bigr\|^{}_1} 
     {\|\vc{x}\|^{}_1} \, .
\end{align*}
  The last term becomes
\[
    \frac{1}{\|\vc{x}\|^{}_1} \bigl\|
      \|\vc{y}\|^{}_1 \vc{x} - \|\vc{x}\|^{}_1 \vc{y} \bigr\|^{}_1 
    \; = \; \frac{1}{\|\vc{x}\|^{}_1} \bigl\| 
      (\|\vc{y}\|^{}_1 - \|\vc{x}\|^{}_1)\ts \vc{x} +
      \|\vc{x}\|^{}_1 (\vc{x} - \vc{y}) \bigr\|^{}_1
    \le 2\ts \|\vc{x}-\vc{y}\|^{}_1 \es,
\]
  from which $\|\cR(\vc{x}) - \cR(\vc{y})\|^{}_1 \le 3 \ts \|\vc{x} -
  \vc{y}\|^{}_1$ follows for $\vc{x}$, $\vc{y} \in \Ll{1}$.  If
  $\vc{x}$, $\vc{y} \in \cM_r$, one has 
  $\|\vc{x}\|^{}_{1} = \|\vc{y}\|^{}_{1}$
  and  the above calculation simplifies to
  $\|\cR(\vc{x}) - \cR(\vc{y})\|^{}_1 
  \le 2\ts \|\vc{x} - \vc{y}\|^{}_1$.  
\hfill \sq

\smallskip \noindent
In continuous time, this is a sufficient condition for the existence
of a unique solution of the initial value problem
\eqref{urconttime}, compare \cite[Thms.~7.6 and 10.3]{Ama}.

\smallskip
It is instructive to generalise the notion of reversibility.
We call a probability vector $\vc{p} \in \cM_1^+$
\emph{reversible} for a recombinator $\cR$ of the form
\eqref{reco} if, for all $i, j, k, \ell \ge 0$,
\begin{equation} \label{urrev}
    T_{ij,k\ell}^{} \, p_k^{} \, p_\ell^{}  \; = \; 
    T_{k\ell,ij}^{} \, p_i^{} \, p_j^{} \, .
\end{equation}
Though this set of equations for detailed balance is much more
restrictive than its linear counterpart in Eq.~\eqref{markov-rev}, the
relevance of this concept is evident from the following property.
\smallskip
\begin{lemma} \label{rev-to-fix}
  If\/ $\vc{p} \in \cM_1^+$ is reversible for\/ $\cR$, it is also a
  fixed point of\/ $\cR$.
\end{lemma}
\Proof  Assume $\vc{p}$ to be reversible for $\cR$.  Then, by \eqref{sumt},
  \begin{displaymath}
    \cR(\vc{p})_i^{} \; = 
    \sum_{j,k,\ell\ge0} \, T_{ij,k\ell}^{} \, p_k^{} \, p_\ell^{} 
    \; =
    \sum_{j,k,\ell\ge0} \, T_{k\ell,ij}^{} \, p_i^{} \, p_j^{} = 
    p_i^{} \sum_{j\ge0} p_j^{} \; = \; p_i^{} \, ,
  \end{displaymath}
for all $i\in\Nnull$, which shows the claim. \hfill \sq

\smallskip \noindent
Returning to the original question of the existence of fixed points,
we now recall the following facts, compare \cite{Bill,Shir} for
details and proofs, which are needed for some general statements 
in the fixed point discussion.
\smallskip
\begin{prop}{\rm \cite[Cor.\ to Thm.~V.1.5]{Yos}}
  \label{vaguenorm}
  Assume the sequence\/ $\bigl(\vc{p}^{(n)}\bigr)$ in\/ $\cM_1^+$ to
  converge in the weak-$*$ topology $($i.e., pointwise, or
  vaguely\/$)$ to some\/ $\vc{p} \in \cM_1^+$, i.e.,
  \begin{displaymath}
    \lim_{n\to\infty} p^{(n)}_k = p_k^{}
    \quad\text{for all\/ $k\in\Nnull$} \es,
    \qquad\text{with\/ $p_k^{}\ge0$ and\/ 
      $\textstyle\sum_{k\ge0} p_k^{} = 1$} \es.
  \end{displaymath}
  Then, it also converges weakly $($in the probabilistic sense\/$)$ and
  in total variation, i.e., $\lim_{n\to\infty} \|\vc{p}^{(n)} -
  \vc{p}\|^{}_1 = 0$.\hspace*{\fill}    \hfill \sq
\end{prop}

\smallskip Recall from \cite{Bill} that a set of measures $\cM \subset
\cM_1^+$ is called \emph{tight} when, for every $\eps>0$, there is an
$m\in\Nnull$ such that $\sum_{k \ge m} p_k^{} < \eps$, simultaneously
for all $\vc{p} \in \cM$.  This is a uniformity condition which serves
as a condition for the compactness needed later on.

\smallskip
\begin{prop}
  \label{prop:urgfixed}
  Assume that the recombinator\/ $\cR$ from\/ $\eqref{reco}$ satisfies
  the normalisation\/ $\eqref{sumt}$ and possesses a convex, weak-$*$
  closed invariant set\/ $\cM \subset \cM_1^+$, i.e., $\cR(\cM)
  \subset \cM$, that is tight.  Then, $\cR$ has a fixed point in\/
  $\cM$.
\end{prop}
\Proof Prohorov's theorem \cite[Thm.~III.2.1]{Shir} states that
  tightness and relative compactness in the weak-$*$ topology are
  equivalent (see also \cite[Chs.~1.1 and 1.5]{Bill}).  In our case,
  $\cM$ is tight and weak-$*$ closed, therefore, due to
  Proposition~\ref{vaguenorm}, norm compact.  Further, $\cM$ is 
  convex by assumption, and 
  $\cR$ is (norm) continuous by Proposition~\ref{recolipschitz}.
  Thus, the claim follows from the Leray--Schauder--Tychonov fixed
  point theorem \cite[Thm.~V.19]{ReSi}. 
\hfill \sq

\smallskip \noindent For several explicit models, we shall see that
such compact invariant subsets indeed exist. On the other hand, once
again due to the infinite-dimensional nature of the dynamical system,
their identification and use for the various proofs is essential.

\section*{5. Takahata's model.}

An early and now classic example was given by Takahata \cite{Tak}.  In
our terminology, he used a recombinator based upon the transition
probabilities
\begin{equation} \label{taka-def}
   T_{ij , k\ell} \; := \; \frac{1}{k+\ell+1} \,
   \delta_{i+j,k+\ell}\, ,
\end{equation}
for $i,j,k,\ell \in \Nnull$. Observing $\mathrm{card} \{(i,j)\mid i,j
\in\Nnull , \, i+j = k+\ell \} = k+\ell+1$, it is clear that $T$ just
describes a recombination with uniform distribution of the copy number
pairs $(k,\ell)$ on each (finite) block of possibilities with $k+\ell$
fixed. One can also check, via Eq.~\eqref{urcons} and
Lemma~\ref{mean-preservation}, that the mean $\m$ is preserved. On the
basis of Eq.~\eqref{taka-def}, the action of the recombinator from
Eq.~\eqref{reco} on probability vectors $\vc{p}$ simplifies to
\begin{equation} \label{taka-reco}
   \cR (\vc{p})_i \; = \; 
   \sum_{\substack{k,\ell\geq 0 \\[0.5mm] k+\ell\geq i}}
   \frac{1}{k+\ell +1}\, p^{}_k \ts p^{}_{\ell} \, .
\end{equation}
Though this model is mathematically rather transparent, it lacks a
good intuitive justification on the level of the biological processes.
Nevertheless, its properties seem to be in acceptable agreement with
at least some of the observations, compare \cite{KDSA,BWZBS}, though
other results, as those shown in \cite{CS}, indicate that also other
types of equilibria appear in experiment.

\smallskip
\begin{prop} \label{taka-prop}
  The probability vector $\bs{p}$ defined by
\[
    p^{}_{n} \; = \; \frac{1}{\m + 1} \,
    \left( \frac{\m}{\m + 1} \right)^n \, , \quad n\in\Nnull \, ,
\]
  is a reversible equilibrium with mean $\m$ for the dynamics
  based on $T$ of\/ $\eqref{taka-def}$. 
\end{prop}
\Proof Using standard identities with geometric series and their
derivatives, it is easy to check that $\bs{p}$ indeed defines a
probability vector on $\Nnull$ with mean $\m$. Detailed balance
follows from a simple calculation,
\[
   \begin{split}
   T^{}_{ij,k\ell}\,p^{}_{k} p^{}_{\ell} & \, = \,
   \frac{\delta_{i+j,k+\ell}}{k+\ell+1}
   \left(\frac{1}{\m +1}\right)^2 
   \left(\frac{\m}{\m +1}\right)^{k+\ell}  \\
   & \, = \, \ts \frac{\delta_{k+\ell,i+j}}{i+j+1}
   \left(\frac{1}{\m +1}\right)^2 \left(\frac{\m}{\m +1}\right)^{i+j}
   \, = \, T^{}_{k\ell,ij}\,p^{}_{i} p^{}_{j} \, ,
   \end{split}
\]
thus completing the claim by means of Lemma~\ref{rev-to-fix}. \hfill \sq

\smallskip These equilibria are geometric distributions as also
discussed above in the Markov context. However, in view of some
experimental findings reported in \cite{CS} and further arguments put
forward in \cite{SA}, one would like to see an initial rise, and
perhaps also a maximum in the vicinity of $n \approx \m$. One should
note that measurements often skip the entries for small copy numbers
(which seem to be rather unreliable), so that a graph with a power law
decay need not indicate the absence of some (weak) form of a maximum.
As the methods for the further analysis of Takahata's model are
similar to what we need later on for alternative models, we first
continue to investigate Takahata's model.

\smallskip
\begin{thm} \label{thm-internal}
  If the initial condition, with mean $\m$, satisfies a certain
  tightness condition $(\limsup_{k\to\infty}\sqrt[k]{p^{}_{k} (0)} <
  1)$, the dynamics, both in discrete and in continuous time,
  converges to the equilibrium vector $\bs{p}$ from
  Proposition~$\ref{taka-prop}$, with $\lim_{t\to\infty} \| 
  \bs{p} (t) - \bs{p} \|_{1} = 0$.
\end{thm}

\smallskip \noindent 
The proof of this theorem, quite appropriately for the present
context, uses an approach via generating functions and then relies on
Banach's contraction principle. It requires several preparatory steps.

\smallskip
Let $\alpha$ and $\delta$ be fixed, with $0<\alpha\le\delta<\infty$,
and consider the space
\begin{equation}  \label{define-X}
  X_{\alpha,\delta} \; := \; \{ \vc{a} = (a^{}_{k})^{}_{k\in\Nnull} \mid
     a^{}_{0}=1\, , \; a^{}_{1}=\alpha\, , \mbox{ and }
     0\le a^{}_{k}\le\delta^k \mbox{ for all } k\ge 2\}\, .
\end{equation}    
If equipped with the metric
\begin{equation} \label{metric-two}
   d(\vc{a},\vc{b}) \; = \; \sum_{k\ge 0} d^{}_{k}\,
   \lvert a^{}_{k} - b^{}_{k}\rvert\, ,
\end{equation}
where $d_k = (\gamma/\delta)^k$ for some $0<\gamma<\tfrac{1}{3}$,
the space $X_{\alpha.\delta}$ is compact \cite[Prop.~5]{RB}.   

Let us define a new vector, $b(\vc{p})$, for suitable $\vc{p}$, by
\begin{equation} \label{def-bmap}
   b(\vc{p})^{}_{k} \; := \; \sum_{\ell\geq k}
   \binom{\ell}{ k}\, p^{}_{\ell} \, ,
\end{equation}
which is certainly well-defined for all $\vc{p}$ with
$\limsup_{k\to\infty} \sqrt[k]{p^{}_{k}} < 1$, by an application of
\cite[Prop.~6]{RB}. This proposition also clarifies the connection
with the space $X_{\alpha,\delta}$ for suitable parameters $\alpha$
and $\delta$. As we shall see, $X_{\alpha,\delta}$ is an example of a
compact, convex space that is invariant under the recombinator
dynamics. It is easy to check that one has $b (\vc{p})^{}_{0} = 1$ and
$b (\vc{p})^{}_{1} = \m$, so that we need $X_{\alpha,\delta}$ with
$\alpha=\m$ and $\delta \ge \m$.

\smallskip   
\begin{lemma} \label{reco-convolve}
  For any $\vc{p}$ with $\limsup_{k\to\infty} \sqrt[k]{p^{}_{k}} < 1$, 
  one has the convolution identity
\[
   b\bigl(\cR (\vc{p})\bigr)_{k} \; = \; \frac{1}{k+1}
   \sum_{m=0}^{k} b(\vc{p})^{}_{m} \ts b(\vc{p})^{}_{k-m} \, .
\]
\end{lemma}
\Proof Let $\vc{p}$ be an arbitrary probability vector with
$\limsup_{k\to\infty} \sqrt[k]{p^{}_{k}} < 1$, so that the mapping
$b$ is well-defined. The left hand side leads to
\[ \begin{split}
   b\bigl(\cR (\vc{p})\bigr)_{k} & \, = \, \sum_{\ell\geq k}
   \binom{\ell}{ k}\,\cR (\vc{p})^{}_{\ell} \, = \,
   \sum_{\ell\geq k} \binom{\ell}{ k} 
   \sum_{\substack{r,s\ge 0 \\[0.5mm] r+s\ge \ell}}
   \frac{p^{}_{r}\ts p^{}_{s}}{r+s+1} \\
   & \, = \, \sum_{\substack{r,s\ge 0 \\[0.5mm] r+s\ge k}}
   \frac{p^{}_{r}\ts p^{}_{s}}{r+s+1} \, \sum_{\ell=k}^{r+s} 
   \binom{\ell}{k} \, = \, \frac{1}{k+1}
   \sum_{\substack{r,s\ge 0 \\[0.5mm] r+s\ge k}}
   \binom{r+s}{k}\, p^{}_{r}\ts p^{}_{s} \, ,
   \end{split}
\]
where a standard identity on binomial coefficients was used
in the last step.

On the other hand, one finds
\[ \begin{split}
   \sum_{m=0}^{k} b(\vc{p})^{}_{m}\ts b(\vc{p})^{}_{k-m}
   & \, = \, \sum_{\substack{m,n\ge 0 \\[0.5mm] m+n=k}}
   \Biggl(\sum_{r\ge m} \binom{r}{m} \, p^{}_{r}\Biggr)
   \Biggl(\sum_{s\ge n} \binom{s}{n} \, p^{}_{s}\Biggr) \\
   & \, = \, \sum_{\substack{m,n\ge 0 \\[0.5mm] m+n=k}}\;
   \sum_{\substack{r,s \ge 0 \\[0.5mm] \, r+s\ge k}} 
   \binom{r}{m} \binom{s}{n}\, p^{}_{r}\ts p^{}_{s} 
   \, = \, \sum_{\substack{r,s \ge 0 \\[0.5mm] 
   r+s\ge k}}  p^{}_{r}\ts p^{}_{s}
   \sum_{\substack{m,n\ge 0 \\[0.5mm] m+n=k}} 
   \binom{r}{m} \binom{s}{n} \\[1mm]
   & \, = \,
   \sum_{\substack{r,s \ge 0 \\ r+s\ge k}} 
   \binom{r+s}{k}\, p^{}_{r}\ts p^{}_{s} \, ,
\end{split}
\]
again using a standard identity, together with the fact that
$\binom{n}{m}=0$ for $m > n$ when $n$ is an integer. A comparison of
the two calculations establishes the claim.  \hfill \sq

\smallskip \noindent
The further relevance of Lemma~\ref{reco-convolve} stems from the 
following property of the generating function of $\vc{p}$, defined by
$\psi(z) = \sum_{\ell\ge 0} p^{}_{\ell} z^{\ell}$.  When rewritten
as a Taylor series around $1$ rather than around $0$, one obtains
\begin{equation}
   \psi(z) \, = \, \sum_{\ell\ge 0} p^{}_{\ell} z^{\ell}
   \, = \, \sum_{k\ge 0}\, \Biggl( \, \sum_{\ell\ge k}
   \binom{\ell}{k}\, p^{}_{\ell} \Biggr) \, (z-1)^k \, = \,
   \sum_{k\ge 0} b(\vc{p})^{}_{k}\, (z-1)^k .
\end{equation}
Under the assumptions on $\vc{p}$, the radius of convergence of
$\psi(z)$ is larger than $1$, so that this calculation is on firm
grounds. Lemma~\ref{reco-convolve} now tells us that we may study
the recombination action on the level of the expansion coefficients.

Let us therefore define the induced recombination operator
$\widetilde{\cR}$ on any space of type $X_{\alpha,\delta}$, with
$\delta\ge\alpha$, by $\widetilde{\cR} (b(\vc{a})) = b(\cR(\vc{a}))$,
which establishes a commuting diagram of the mappings $\cR$ and
$\widetilde{\cR}$ in the obvious way. More precisely, one first
restricts the action of $\cR$ to a suitable subspace of $\cM^{+}_{1}$,
so that the mapping $b$ is well-defined.  If $\vc{p}$ satisfies the
condition of Lemma~\ref{reco-convolve}, so that the radius of
convergence of $\psi$ exceeds $1$, the probability vector $\vc{p}$ is
also completely determined by its moments, compare
\cite[Thm.~II.12.7]{Shir} together with the observation that
$\psi(e^{it})$ is the (convergent) moment generating function of
$\vc{p}$. As all moments, in turn, are specified by the entries of
$b(\vc{p})$, the latter uniquely determines $\vc{p}$ in this
situation.

It is easy to check that the vector $(1,\alpha,\alpha^2,\ldots )$ is a
fixed point of $\widetilde{\cR}$ in $X_{\alpha,\delta}$, for any
$\delta\ge\alpha$.  Choosing $\alpha=\m$, this vector is the image of
the probability vector $\vc{p}$ from Proposition~\ref{taka-prop} under
the mapping $b$.

\smallskip
\begin{prop}
  On $X_{\alpha,\delta}$, the map defined by $\widetilde{\cR}$ is a
  contraction. In particular, it is a globally Lipschitz continuous
  mapping of $X_{\alpha,\delta}$ into itself.
\end{prop}
\Proof Let $\delta\ge \alpha> 0$ be given, as well as arbitrary
$\vc{a},\vc{b}\in X_{\alpha,\delta}$. Clearly, we have
$\widetilde{\cR} (\vc{a})^{}_{0}=1$ and $\widetilde{\cR}
(\vc{a})^{}_{1}=\alpha$.  For $k\ge 2$, one finds $\widetilde{\cR}
(\vc{a})^{}_{k}=\frac{1}{k+1}\sum_{\ell=0}^{k} a^{}_{\ell}\ts
a^{}_{k-\ell} \le \delta^k$. This proves that $\widetilde{\cR}$ maps
$X_{\alpha,\delta}$ into itself.

The space $X_{\alpha,\delta}$ is equipped with the metric $d$ from
\eqref{metric-two}. Since, due to $\vc{b}\in X_{\alpha,\delta}$, also
$\widetilde{\cR} (\vc{b})^{}_{0} = 1$ and $\widetilde{\cR}
(\vc{b})^{}_{1} =\alpha$, the contraction estimate reads as follows.
\[  \begin{split}
   d(\widetilde{\cR} (\vc{a}), \widetilde{\cR} (\vc{b})) & \, = \, 
   \sum_{k\ge 2} \frac{d_k}{k \! + \! 1}
   \,\Bigl\lvert \sum_{\ell=0}^{k} (a^{}_{\ell}\ts a^{}_{k-\ell}
   \!\ts - b^{}_{\ell}\ts b^{}_{k-\ell} ) \Bigr\rvert 
   \, = \ts \sum_{k\ge 2} \frac{d_k}{k \! + \! 1}
   \,\Bigl\lvert \sum_{\ell=0}^{k} (a^{}_{\ell} \! - b^{}_{\ell})
   (a^{}_{k-\ell}\!\ts + b^{}_{k-\ell} ) \Bigr\rvert \\
   & \,\le \, \sum_{k\ge 2} \frac{2\, d_k}{k+1}
   \sum_{\ell=2}^{k} \delta^{k-\ell} \, \lvert
   a^{}_{\ell} - b^{}_{\ell}\rvert 
   \, = \, \sum_{\ell\ge 2} d^{}_{\ell}\,
   \lvert a^{}_{\ell} - b^{}_{\ell}\rvert
   \sum_{k\ge \ell} \frac{2}{k+1}\, \delta^{k-\ell}\,
   \frac{d_k}{d_{\ell}} \, .
\end{split}
\]
With the choice $d_k = (\gamma/\delta)^k$, where we had
$\gamma < \frac{1}{3}$,  we can now find, for $\ell\ge 2$,
an upper bound for the inner sum,
\[
   \sum_{k\ge \ell} \frac{2}{k+1}\, \delta^{k-\ell}\,
   \frac{d_k}{d_{\ell}} \, \le \, \frac{2}{3}
   \sum_{k\ge \ell} \gamma^{k-\ell} \, = \,
   \frac{2}{3-3\gamma} \, =: \, C \, < \, 1\, ,
\]
which, together with the previous calculation, proves the contraction
property,
\[
   d(\widetilde{\cR} (\vc{a}), \widetilde{\cR} (\vc{b}))
   \, \le \, C\, d(\vc{a},\vc{b})\, ,
\]
with contraction constant $C < 1$. Clearly, this also means that
$\widetilde{\cR}$ is globally Lipschitz continuous.  \hfill \sq

\smallskip This shows that, in discrete time, we have exponentially
fast convergence of the sequence $( \widetilde{\cR} ^n (\vc{a}))_{n\ge
  1}$, with $\vc{a}\in X_{\alpha,\delta}$, to a unique fixed point in
$X_{\alpha,\delta}$. It is specified by the mean copy number $\m$ of
the probability vector $\vc{p}$ that underlies $\vc{a}=b(\vc{p})$, via
$\alpha = \m$, see above. Clearly, this fixed point (in
$X_{\alpha,\delta}$) is the image (under $b$) of the equilibrium
vector $\vc{p}\in\cM^{+}_{1}$ calculated earlier in
Proposition~\ref{taka-prop}, as the mapping $b$ is invertible in this
situation.  The claim of Theorem~\ref{thm-internal} for discrete time
is now clear, with exponentially fast convergence to the equilibrium,
from any initial condition as specified there.

{}For the slightly more involved treatment of the continuous
time case, we refer to \cite{RB}. It is based on the construction of a
Lyapunov function, similar to that of \cite[Prop.~13]{RB}.

\section*{6. Internal crossover.}

Another rather obvious model is based on the assumption that the
shorter of the two sequences (or stretches) can align with any
connected block of the longer sequence, but without any overhang.
This situation has been coined internal unequal crossover, or
\emph{internal crossover} for short. Here, restricting to probability
measures on $\Nnull$, the recombinator \eqref{reco} simplifies to
\begin{equation} \label{reco-internal}
   \cR_{0} (\vc{p})_i \; = \; \sum_{\substack{k,\ell\ge 0 \\
   k \land \ell \le i \le k \lor \ell}} 
   \frac{p^{}_{k}\ts p^{}_{\ell}}{1+\lvert k-\ell\rvert} \, ,
\end{equation}
where $k\land\ell$ ($k\lor\ell$) stands for the minimum (the maximum)
of $k$ and $\ell$, see \cite{SA,RB,Oli} for details on this model.
We choose the notation $\cR_0$ for reasons that will become clear
later on.

In our search for fixed points, it is again useful to look for 
probability vectors that are reversible for $\cR_0$.
Since both forward and backward transition
probabilities are simultaneously non-zero only when $\{i,j\} =
\{k,\ell\} \subset \{n,n+1\}$ for some $n$, the components $p_k^{}$
may only be positive on this small set as well.  By the following
proposition, this indeed characterises all fixed points of this case.

\smallskip
\begin{prop} \label{fixed-internal} A probability measure $\vc{p} \in
  \cM_1^+$ is a fixed point of\/ $\cR_0$ if and only if its mean copy
  number\/ $\m = \sum_{k\ge0} k\,p_k^{}$ is finite, together with
  $p_{\lfloor \m \rfloor}^{} = \lfloor \m \rfloor + 1 - \m$,
  \linebreak $p_{\lceil \m \rceil}^{} = \m + 1 - \lceil \m \rceil$,
  and $p_k^{} = 0$ for all other $k$.  This includes the case that\/
  $\m$ is a non-negative integer, where $p_{\lfloor \m \rfloor}^{} =
  p_{\lceil \m \rceil}^{} = p_{\m}^{} = 1$.
\end{prop}
\Proof The `if' follows easily by insertion into \eqref{urrev} and
Lemma~\ref{rev-to-fix}.  For the `only if' part, let $i$ denote the
smallest integer such that $p_i^{} > 0$.  Then,
\[
    \cR (\vc{p})_i^{} \, = \, 
    p_i^2 + 2 p_i^{} \sum_{\ell\ge1} \frac{p_{i+\ell}^{}}{1+\ell} 
    \, = \,
    p_i^{} \left( p_i^{} + p_{i+1}^{} + 
      \sum_{\ell\ge2} \frac{2}{\ell+1} p_{i+\ell}^{} \right) 
    \, \le \,   p_i^{} \, ,
\]
where the last step follows since $\frac{2}{\ell+1} < 1$ in the last
sum. One has equality precisely when $p_k^{} = 0$ for all $k \ge i+2$.
This implies $\m < \infty$ and the uniqueness of $\vc{p}$ (given $\m$)
with the non-zero frequencies as claimed.
\hfill \sq

\smallskip
In this case, one may select a compact subset within the probability 
vectors by demanding the existence of the centred $r$-th
moment, for some fixed $r>1$. More precisely, with
\[
    \mu^{}_s (\vc{p}) \; := \; \sum_{\ell\ge 0}\; 
   \lvert \ell -\m\rvert^s\ts p^{}_{\ell}\, , 
\]
one considers the set
\begin{equation} \label{moment-set}
  \cM^{+}_{1,\m,C} \; := \;
   \{\vc{p}\in\cM^{+}_{1} \mid \mbox{$\sum_k k\ts p^{}_{k} = \m$  and }
    \mu^{}_r (\vc{p}) \le C \}
\end{equation}   
for an arbitrary, but fixed $C < \infty$, equipped with our usual
metric as introduced before in \eqref{urgpmetric}. This gives
a compact and convex space \cite[Lemma~2]{RB}.  Moreover, one has

\smallskip
\begin{lemma}  \label{lyapunov-1}
   Let $r>1$ be fixed and consider the space $\cM^{+}_{1,\m,C}$ of\/
   $\eqref{moment-set}$. Then, both $\mu^{}_1$ and $\mu^{}_r$ satisfy
\[
    \mu^{}_s (\cR_0 (\vc{p})) \; \le \; \mu^{}_s (\vc{p}) \, ,
\]
   with equality if and only if $\vc{p}$ is a fixed point of $\cR_0$.
   
   Moreover, $\mu^{}_1 \! : \; \cM^{+}_{1,\m,C} \longrightarrow
   \RR_{\ge 0}$ is continuous and defines a Lyapunov function for the
   dynamics in continuous time.
\end{lemma}       
\Proof To show the first claim, consider
\begin{equation}    \label{moment-estimates}
    \begin{aligned}
      \mu^{}_s(\cR_0(\vc{p})) &\; = \;
      \sum_{i\ge0} \!
        \sum_{\substack{k,\ell\ge0 \\[0.4\baselineskip] 
            k\wedge\ell \le i \le k\vee\ell}} \!
        \frac{|i-\m|^s}{1+|k-\ell|} \, p_k^{} \, p_\ell^{} \\
      &\; = \; \sum_{k,\ell\ge0} \frac{p_k^{} \, p_\ell^{}}{1+|k-\ell|}
        \, \frac12 \sum_{i=k\wedge\ell}^{k\vee\ell}
        (|i-\m|^s + |k+\ell-i-\m|^s) \, .
    \end{aligned}
  \end{equation}
  For notational convenience, let $j = k+\ell-i$.  We now show
  \begin{equation}   \label{help-one}
    |i-\m|^s + |k+\ell-i-\m|^s \; \le \; |k-\m|^s + |\ell-\m|^s \, .
  \end{equation}
  If $\{k,\ell\} = \{i,j\}$, then \eqref{help-one} holds with
  equality.  Otherwise, assume without loss of generality that $k < i
  \le j < \ell$.  If $\m \le k$ or $\m \ge \ell$, we have equality for
  $s=1$, but a strict inequality for $s=r$ due to the convexity of $x
  \mapsto x^r$.  (For $s=1$, this describes the fact that a
  recombination event between two sequences that are both longer or
  both shorter than the mean does not change their averaged distance
  to the mean copy number.)  In the remaining cases, the inequality is
  strict as well.  Hence, $\mu^{}_s(\cR_0(\vc{p})) \le \mu^{}_s
  (\vc{p})$ with equality if and only if $\vc{p}$ is a fixed point of
  $\cR_0$, since otherwise the sum in \eqref{moment-estimates}
  contains at least one term for which \eqref{help-one} holds as a
  strict inequality.

  To see that $\mu^{}_1$ is continuous, consider a convergent sequence
  $(\vc{p}^{(n)})$ in $\cM^{+}_{1,\m,C}$ and the random variables
  ${H}^{(n)} = |K^{(n)}-\m|$, where the $K^{(n)}$ are independent
  $\Nnull$-valued random variables with laws $\vc{p}^{(n)}$.  Due to
  the structure of $\cM^{+}_{1,\m,C}$, the random variables $H^{(n)}$
  are uniformly integrable, which implies the convergence of the
  corresponding expectation values by \cite[Thm.~25.12]{B1}.  This, in
  turn, is nothing but the continuity of $\mu^{}_1$.  Since
  $\mu^{}_1(\vc{p})$ is linear in $\vc{p}$ and thus infinitely
  differentiable, so is the solution $\vc{p}(t)$ for every initial
  condition $\vc{p}_0 \in \cM^{+}_{1,\m,C}$, compare \cite[Thm.~9.5
  and Remark~9.6(b)]{Ama}.  Therefore, we have
\[
    \dot{\mu^{}}_1 (\vc{p}_0) \, = \,
    \liminf_{t\to0^+} \frac{\mu^{}_1(\vc{p}(t)) - 
    \mu^{}_1(\vc{p}_0)}{t} \,  = \,
    \mu^{}_1 (\cR^{}_0(\vc{p}_0)) - \mu^{}_1 (\vc{p}_0) \, \le \, 0 \, ,
\]
again with equality if and only if $\vc{p}_0$ is a fixed point.
Thus, $\mu^{}_1$ is a Lyapunov function as claimed.
\hfill \sq

\smallskip
Finally, this gives the following convergence result, the proof of which 
is given in \cite{RB} and not repeated here.
\smallskip
\begin{thm}  \label{internal-reco-thm}
 Assume that, for the initial condition $\vc{p}(0)$ and fixed
  $r>1$, the $r$-th moment exists, $\mu^{}_r (\vc{p}) <
  \infty$.  Then, $\m = \sum_{\ell} \ell \ts p^{}_{\ell}$ is finite and,
  both in discrete and in continuous time, $\lim_{t\to\infty}
  \|\vc{p}(t) - \vc{p}\|^{}_1 = 0$ with the appropriate fixed point
  $\vc{p}$ from Proposition~$\ref{fixed-internal}$.
  \hfill \sq
\end{thm} 

\smallskip Let us mention that, for $q=0$, the recombinator can be
expressed in terms of explicit frequencies $\pi_{k,\ell}$ of fragment
pairs before concatenation (with copy numbers $k$ and $\ell$) as
$\cR_0(\vc{p})_i^{} = \sum_{j=0}^i \pi_{j,i-j}$.  It is as yet an open
question whether this can be used to simplify the above treatment.

\section*{7. Random crossover.}

This model deviates from the previous one in that it admits arbitrary
overhangs, up to the case where, after the crossover, one sequence got
it all while the other lost everything. The possible alignments for
any pair are supposed to be equally likely, so that the recombinator
\eqref{reco}, again restricted to the probability measures, now reads
\begin{equation} \label{reco-random}
    \cR_1(\vc{p})_i^{} \; = \;
  \sum_{\substack{k,\ell\ge 0 \\[0.5mm] \; k+\ell \ge i}}
    \frac{1+\min\{k,\ell,i,k+\ell-i\}}{(k+1)(\ell+1)} 
    \, p_k^{} \, p_\ell^{} \, .
\end{equation}
As for our previous two examples, using Lemma~\ref{rev-to-fix} once
again, the reversibility condition,
\[
  \frac{p_{k\vphantom{j}}^{}}{k+1} \frac{p_{\ell\vphantom{j}}^{}}{\ell+1} 
  \; = \; \frac{p_{i\vphantom{j}}^{}}{i+1} \frac{p_j^{}}{j+1} \, ,
  \quad\text{for all } k+\ell = i+j \, ,
\]
leads to an expression for fixed points. In fact, these relations have
$p_k^{} = C (k+1) x^k$ as a solution, with appropriate parameter $x$
and normalisation constant $C$. Again, it turns out that all fixed
points are given this way, as was originally noticed (in a different
way) in \cite[Thm.~A.2]{SA}.

\smallskip
\begin{prop} \label{random-fixed}
  Every fixed point $\vc{p} \in \cM_1^+$ of\/ $\cR_1$ has finite mean
  $\m = \sum_k k \ts p^{}_{k}$, and is uniquely specified by the value of
  $\m$. Explicitly, one has
\[
    p_k^{} \; = \; \left(\frac{2}{\m+2}\right)^2 (k+1) 
       \left(\frac{\m}{\m+2}\right)^k  
\]
  with $k \in \Nnull $.       \hfill \sq
\end{prop}

\smallskip  \noindent
One can verify this in several ways, one being a direct calculation
via induction. Interestingly, this equilibrium is the convolution
of two geometric distributions (of equal mean $\m /2$), and hence also
of geometric type according to our terminology (which follows that
of \cite{SA}). It might be interesting to explore this observation a
little further in the future.

At this point, one can define, very much in analogy to the situation
in Takahata's model above, an induced recombinator, $\widetilde{\cR}_{1}$,
acting once more on spaces of the form $X_{\alpha,\delta}$. It is
given as 
\[
     \widetilde{\cR}_1 (\vc{p}) \; = \; a \bigl(\cR_1 (\vc{p})\bigr)
\]
where the mapping $a$ is defined by
\[
     a (\vc{p})^{}_k \; = \; \frac{1}{k+1} \sum_{\ell\ge k}
     \binom{\ell}{k} p^{}_{\ell}  \; = \; \frac{1}{k+1}\,
     b (\vc{p})^{}_k \, .
\]
It is thus closely related to our above mapping $b$.

The main result on this model, proved in detail in \cite{RB,Oli},
reads as follows.
\smallskip
\begin{thm} \label{random-convergence} 
  Assume that\/ $\,\limsup_{k\to\infty} \sqrt[k]{p_k^{}(0)} < 1$.
  Then, both in discrete and in continuous time,\/ $\lim_{t\to\infty}
  \|\vc{p}(t) - \vc{p}\|^{}_1 = 0$, where $\vc{p}$ is the
  corresponding fixed point according to
  Proposition~$\ref{random-fixed}$. 
\end{thm}
\Proof The proof is very similar to the one used above for the
Takahata model, and employs once again Banach's contraction principle
for the induced action of $\widetilde{\cR}_{1}$ on
$X_{\alpha,\delta}$. Since all details have been given in \cite{RB},
we omit them here.  \hfill \sq

\smallskip The fixed points of Proposition~\ref{random-fixed} are of
the expected geometric type, and are perhaps more realistic than those
of the Takahata model, at least for cases where a maximum is present
in the repeat distribution. However, one should note that the
experimental situation is not completely convincing at present, so
that it seems advantageous to have a versatile model class at hand.

\section*{8. An interpolation.}

When considering the recombinators $\cR_0$ and $\cR_1$ in comparison,
one would like to find further models that share properties of both of
them, or interpolate between them in a suitable way. In particular,
$\cR_0$ is unrealistic due to the complete confinement of the shorter
bit within the range of the longer one, while $\cR_1$ poses no
restriction at all for any kind of overhang. One such interpolation
was initially investigated in discrete time by Atteson and Shpak in
\cite{SA}, based on preceding work by Ohta \cite{Oht} and Walsh
\cite{Wsh}, see also \cite{RB,Oli} for more. The interpolation employs
a penalty function idea for overhangs of the shorter sequence, and
leads (in the above language) to a recombinator $\cR_q$ with $0\le
q\le 1$. The latter is based upon the transition probabilities
\begin{equation} \label{interpol-one}
   T^{(q)}_{ij,k\ell} = 
  C^{(q)}_{k\ell} \, \delta_{i+j,k+\ell} \, (1+\min\{k,\ell,i,j\}) \,
    q^{0 \vee (k\wedge\ell - i \wedge j)} \, ,
\end{equation}
where $k \vee \ell := \max\{k,\ell\}$, $k \wedge \ell :=
\min\{k,\ell\}$, and $0^0 = 1$.  The normalisation constants
$C^{(q)}_{k\ell}$ are chosen such that \eqref{sumt} holds, i.e., 
$\sum_{i,j\ge0} T^{(q)}_{ij,k\ell} = 1$. These constants are symmetric in 
$k$ and $\ell$ and read explicitly
\[
  C^{(q)}_{k\ell} =
  \frac{(1-q)^2}{
    (k\wedge\ell+1)(|k-\ell|+1)(1-q)^2 +
    2q(k\wedge\ell - (k\wedge\ell+1)q + q^{k\wedge\ell+1})} \, .
\]
Note further that the total number of units is indeed conserved in
each event and that the process is symmetric within both pairs.  Hence
\eqref{urcons} is satisfied.

Unfortunately, the situation with the fixed points is a lot more
complicated due to the following result.

\smallskip
\begin{prop} \label{no-rev}
  For parameter values $q\in (0,1)$, any fixed point $\vc{p}
  \in \cM_1^+$ of the recombinator $\cR_q$, given by $\eqref{reco}$ and
  $\eqref{interpol-one}$, satisfies $p_k^{}>0$ for all $k\ge0$ $($unless it
    is the trivial fixed point $\vc{p} = (1,0,0,\ldots)$ we
    excluded\/$)$.  None of these extra fixed points is reversible.
\end{prop}
\Proof  Let a non-trivial fixed point $\vc{p}$ be given and choose any $n>0$
  with $p_n^{}>0$.  Observe that $T^{(q)}_{n+1 \; n-1, nn} > 0$ for
  $0 < q < 1$ and hence
  \begin{displaymath}
    p_{n\pm1}^{} \, = \, \cR_q(\vc{p})_{n\pm1}^{} 
    \, = \,  \sum_{j,k,\ell\ge0} T^{(q)}_{n\pm1 \; j,k\ell} \, 
      p_k^{} \, p_\ell^{}
    \, \ge \, T^{(q)}_{n+1 \; n-1, nn} \, p_n^{} \, p_n^{} > 0 \, .
  \end{displaymath}
  The first statement now follows by induction.  For the second statement,
  evaluate the reversibility condition \eqref{urrev} for all
  combinations of $i$, $j$, $k$, $\ell$ with $i+j = k+\ell \le 4$.
  This leads to four independent equations.  Three of them can be
  transformed to the recursion
 \[
    p_k^{} \, = \, \frac{(k+1)q}{2(k-1)+2q} 
    \frac{p_1^{}}{p_0^{}} \, p_{k-1}^{} \, ,
    \qquad k \in \{2,3,4\} \, ,
 \]
  from which one derives explicit equations for all $p_k^{}$ with $k
  \in \{2,3,4\}$ in terms of $p_0^{}$ and $p_1^{}$.  Inserting the one
  for $p_2^{}$ into the remaining equation yields another equation for
  $p_4^{}$ in terms of $p_0^{}$ and $p_1^{}$, which contradicts the
  first equation for all $q \in (0,1)$, as is easily verified.
\hfill \sq

\smallskip
Nevertheless, the dynamics is well defined, and respects the compact
subsets defined above in forward time, compare \cite[Thm.~4]{RB}.
Based upon the analysis in \cite{Oli,RB}, and further numerical work
on the fixed points, it is plausible that, given the mean copy number 
$\m$, never more than one fixed point for $\cR_q$ exists.  
Due to the global convergence
results at $q=0$ and $q=1$, any non-uniqueness in the vicinity of
these parameter values could only come from a bifurcation, not from an
independent source.  Numerical investigations indicate that no
bifurcation is present, but this needs to be analysed further.

Moreover, the Lipschitz constant for the corresponding induced
recombinator $\widetilde{\cR}_q$ can be expected to be continuous in
the parameter $q$, hence to remain strictly less than $1$ on the sets
$X_{\alpha,\delta}$ in the neighbourhood of $q=1$.  So, at least
locally, the contraction property should be preserved. For further
progress, it seems advantageous \cite{Hof} to use a rather different
approach based on the analysis of similar problems in evolutionary
game theory.  Here, one would aim to establish a slightly weaker type
of convergence result for all $0 < q < 1$, and probably even on the
larger compact set $\cM^{+}_{1,\m,C}$ from Eq.~\eqref{moment-set}.

\section*{9. Open problems and outlook.}

The results for the various models presented here show that initial
configurations, subject to some specific conditions that are no
restriction in practice, converge to one of the known fixed
points. These results apply to the deterministic dynamics of the
infinite population limit.

In view of the biological applications, one is also interested in
possible deviations from this picture on the level of large, but
finite, populations, i.e., for the underlying stochastic process,
e.g., a variant of the Moran model with unequal crossover.  In this
model class, however, important deviations seem unlikely, due to the
known convergence results for the infinite population limit, see
\cite{BaHi} for more.

Since the above equilibrium distributions have finite support or are
exponentially small for large copy numbers, one can also expect these
systems to behave very much like ones with only finitely many
types. In this sense, the results are typical, and the more general
setting with probability vectors on $\Nnull$ is adequate. This is also
supported by several simulations \cite{Oli}.

Still, an open question is a more complete understanding of the regime
$q\in (0,1)$ in Section~8. Due to the loss of reversibility of the
fixed points, the analysis becomes rather involved. Preliminary
investigations \cite{Oli} have not given any hint on values of $q$
where convergence fails or where alternative stable fixed points show
up, though this is presently only based on numerical experiments and
perturbative arguments.  It might be advantageous (and perhaps also
more realistic) to search for other ways to interpolate between the
cases $q=0$ and $q=1$, preferably ones that maintain the reversibility
of the equilibria. This question certainly deserves further attention.



\begin{thebibliography}{99}

\bibitem{Ama}
H.~Amann,
\textit{Ordinary Differential Equations},
de Gruyter, Berlin (1990).

\bibitem{BaHi}
E.~Baake and I.~Herms,
\textit{Single-crossover dynamics:\ finite versus infinite 
populations}, Bull.\ Math.\ Biol.\ \textbf{70} (2008) 603--624;
\texttt{arXiv:q-bio/0612024}.

\bibitem{BWZBS}
D.~Bachtrog, S.~Weiss, B.~Zangerl, G.~Brem and
C.~Schl\"otterer,
\textit{Distribution of dinucleotide microsatellites
in the \emph{Drosophila melanogaster} genome},
Mol.\ Biol.\ Evol.\ \textbf{16} (1999) 602--610.

\bibitem{Bank}
C.~Bank,
\textit{Diskrete Rekombinationsdynamik f\"ur repetitive Strukturen},
Diplomarbeit, Univ.\ Bielefeld (2007).

\bibitem{B1}
P.~Billingsley, 
\textit{Probability and Measure},
3rd ed., Wiley, New York (1995).

\bibitem{Bill}
P.~Billingsley, 
\textit{Convergence of Probability Measures},
2nd ed., Wiley, New York (1999).

\bibitem{Bue}
R.~B\"urger,
\textit{The Mathematical Theory of Selection,
Recombination and Mutation},
Wiley, Chichester (2000).

\bibitem{Bus}
C.\thinspace D.~Bustamante,
\textit{Population genetics of molecular evolution},
in \cite{N}, pp.~63--99.

\bibitem{CS}
P.~Calabrese and R.~Sainudiin,
\textit{Models of microsatellite evolution},
in \cite{N}, pp.~289--305.

\bibitem{GL}
D.~Graur and W.-H.~Li, 
\textit{Fundamentals of Molecular Evolution},
2nd ed., Sinauer, Sunderland (2000).

\bibitem{Hof}
J.~Hofbauer, private communication (2003).

\bibitem{K}
B.~Kitchens,
\textit{Symbolic Dynamics --- One-sided, Two-sided and
Countable State Markov Chains}, 
Springer, Berlin (1998). 

\bibitem{KDSA}
S.~Kruglyak, R.\thinspace T.~Durrett, 
M.\thinspace D.~Schug and C.\thinspace F.~Aquadro,
\textit{Equilibrium distributions of microsatellite
repeat length resulting from a balance between slippage
events and point mutations},
Proc.\ Natl.\ Acad.\ Sci.\ USA \textbf{95}
(1998) 10774--10778.

\bibitem{L1}
M.~Lachowicz, 
\textit{General population systems. Macroscopic limit of
a class of stochastic semigroups},
J.\ Math.\ Anal.\ Appl.~\textbf{307} (2005) 585--605.

\bibitem{L2}
M.~Lachowicz, 
\textit{Micro and meso scales of description corresponding
to a model of tissue invasion by solid tumors},
Math.\ Models Meth.\ Appl.\ Sci.~\textbf{15} (2005) 1667--1683.

\bibitem{N}
R.~Nielsen (ed),
\textit{Statistical Methods in Molecular Evolution},
Springer, New York (2005).

\bibitem{Oht}
T.~Ohta,  
\textit{On the evolution of multigene families},
Theor.\ Pop.\ Biol.~\textbf{23} (1983), 216--240.

\bibitem{Oli} 
O.~Redner,
\textit{Models for Mutation, Selection, and Recombination
in Infinite Populations},
Dissertation, Univ.\ Greifswald;
Shaker, Aachen (2003).

\bibitem{RB}
O.~Redner and M.~Baake,
\textit{Unequal crossover dynamics in discrete and
continuous time},
J.\ Math.\ Biol.~\textbf{49} (2004), 201--226;
\texttt{arXiv:math.DS/0402351}.

\bibitem{ReSi}
M.~Reed and B.~Simon, 
\textit{Methods of Modern Mathematical Physics {I}:\ 
Functional Analysis}, Academic Press, San Diego (1980).

\bibitem{S}
E.~Seneta,
\textit{Non-negative Matrices and Markov Chains},
rev.\ printing, Springer, New York (2006).

\bibitem{Shir}
A.\thinspace N.~Shiryaev,
\textit{Probability},
2nd ed., Springer, New York (1996).

\bibitem{SA}
M.~Shpak and K.~Atteson, 
\textit{A survey of unequal crossover systems and 
their mathematical properties},
Bull.\ Math.\ Biol.~\textbf{64} (2002), 703--746.

\bibitem{Tak}
N.~Takahata,
\textit{A mathematical study on the distribution of the 
number of repeated genes per chromosome},
Genet.\ Res.~\textbf{38} (1981), 97--102.

\bibitem{Wake}
J.~Wakeley,
\textit{Coalescent Theory:\ An Introduction},
Roberts and Company,  Greenwood Village, CO (2008).

\bibitem{Wsh}
J.\thinspace B.~Walsh, 
\textit{Persistence of tandem arrays:\ Implications for 
satellite and simple-sequence {DNA}s},
Genetics \textbf{115} (1987), 553--567.
  
\bibitem{Yos}
K.~Yosida,
\textit{Functional Analysis}, 
6th ed., Springer, Berlin (1980).

\end{thebibliography}
\end{document}